\documentclass[10pt]{amsart}
\usepackage{amsmath, amscd, amsthm, amssymb, amsfonts, verbatim}
\usepackage[all]{xy}

\title{The $A_\infty$ operad and the moduli space of curves}
\author{Kevin Costello}
\address{Department of Mathematics\\
Imperial College London\\ London SW7 2AZ\\ United Kingdom}
\email{k.costello@imperial.ac.uk}
\date{}

\newcommand{\iso}{\cong}
\newcommand{\diag}{\bigtriangleup} 
\newcommand{\C}{\mathbb C}

\newcommand{\Q}{\mbb Q}

\newcommand{\Z}{\mathbb Z}

\newcommand{\into}{\hookrightarrow}

\newcommand{\op}{\operatorname}

\newcommand{\mbb}{\mathbb}

\newcommand{\mc}{\mathcal}

\newcommand{\abs}[1]{\left| #1 \right|}

\newcommand{\R}{\mbb R}

\newcommand{\opmod}{\overline { \mc N}}

\newcommand{\open}{\opmod}

\newcommand{\rrar}{\twoheadrightarrow}

\newcommand{\Graphs}{ {\bf Graphs} }
\newcommand{\Trees} { {\bf Forests} }
\newcommand{\RTrees}{ {\bf RForests} }
\newcommand{\Fin}{ {\bf Fin} }

\newcommand{\Cyc}{ {\bf Cyc} }
\newcommand{\Mod}{ {\bf Mod} }
\newcommand{\Pairs}{ {\bf Pairs} }
\newcommand{\Top}{ {\bf Top} }
\newcommand{\Orb}{ {\bf Orb} }
\newcommand{\Vect}{ {\bf Vect} }

\newcommand{\dgv}{ {\bf dg} }

\newtheorem{theorem}{Theorem}[subsection]
\newtheorem{proposition}[theorem]{Proposition}
\newtheorem{definition}[theorem]{Definition}
\newtheorem{lemma}[theorem]{Lemma}

\newtheorem{corollary}[theorem]{Corollary}

\numberwithin{equation}{subsection}

\begin{document}

\begin{abstract}
The modular envelope of a cyclic operad is the smallest modular operad
containing it.  A modular operad is constructed from moduli spaces of
Riemann surfaces with boundary; this modular operad is shown to be the
modular envelope of the $A_\infty$ cyclic operad. This gives a new
proof of the result of Harer-Mumford-Thurston-Penner-Kontsevich that a
cell complex built from ribbon graphs is homotopy equivalent to the moduli space of curves.
\end{abstract}
\maketitle

\section{Introduction}

In recent years, an interesting relationship has emerged between
the $A_\infty$ (or associative) operad and the moduli spaces of
curves.  Probably the first person to make the connection to the
associative operad was Witten \cite{wit1986, wit1995}, who found
the associativity relation appearing in open string  theory.
Around the same time, Harer-Mumford-Thurston \cite{har1986} and
Penner \cite{pen1987} showed that the cohomology of the moduli
space of curves can be described by a complex built out of ribbon
graphs. It was shown by Kontsevich \cite{kon1994} that ribbon
graphs are closely related to the $ A_\infty$-operad; in
particular he used ribbon graphs to associate to any $
A_\infty$-algebra certain cohomology classes in the moduli space
of curves.

Several aspects of this story are a little unsatisfactory.  For
example, the results of Harer-Mumford-Thurston-Penner and
Kontsevich rely on cell decompositions of the moduli space of
curves, with a cell for each ribbon graph.  However, there are
many such triangulations known, with no canonical choice.  A clear
geometric reason for the description of the cohomology of moduli
space by ribbon graphs seems to be lacking.  Also, ribbon graphs
form a (modular) operad --  ribbon graphs can be glued along
external edges.  It was not clear (at least to me)  what operadic
structure on moduli space corresponds to this structure on ribbon
graphs.

In another direction, it is known that a certain topological
operad constructed from holomorphic discs with marked points on
the boundary is isomorphic to the $A_\infty$ topological operad.
This is the reason for the appearance of $A_\infty$ algebras in
Floer homology and the Fukaya category \cite{fuk_oh_oht_ono2000,
fuk2001, flo1988, flo1989}. The associative operad appears in a
closely related way in the work of Moore and Segal \cite{moo2001,
seg2001}.  They show that an open topological field theory, at all
genera, is given by a (not necessarily commutative) Frobenius
algebra.  This is an analogue of the well-known result that a
closed topological field theory is the same as a commutative
Frobenius algebra.

This note describes a  different point of view on the relationship
between the $A_\infty$ operad and the moduli space of curves,
where all of the above results can be seen naturally.  In
particular, new proofs of the results of Kontsevich and
Harer-Mumford-Thurston-Penner are given.  The main result is that
that the modular operad controlling open topological conformal
field theory, at all genera, is the smallest modular operad
containing the cyclic operad of $A_\infty$ algebras.  One
immediate corollary is that an open topological conformal field
theory is the same as an $A_\infty$ algebra with invariant inner
product. This is a generalisation of the work of Moore and Segal
on open topological
 field theory.  The distinction between topological field theory
and topological conformal field theory is that the former deals
with topological surfaces, while the latter takes account of the
topology of moduli spaces of conformal structures on surfaces.
Another corollary of the main result is the ribbon graph
decomposition of moduli space.

\begin{figure}
$$\begin{xy} 0*\cir<30pt>{}, 0,a(0), **{}, 0+/28pt/*{\bullet},
0+/56pt/="c2"*\cir<30pt>{}, 0,a(120), **{},
 0+/28pt/*{\bullet},
 0,a(240),**{}, 0+/28pt/*{\bullet},
 "c2",a(60), **{},
 "c2"+/28pt/*{\bullet},
 "c2",a(300), **{},
 "c2"+/28pt/*{\bullet},
 0, a(300), **{}, 0+/56pt/="c3"*\cir<30pt>{},
 0+/28pt/*{\bullet}, "c3", a(60), **{},
 "c3"+/28pt/*{\bullet}, "c3", a(270), **{},
 "c3"+/28pt/*{\bullet}
\end{xy}$$
\caption{A point in $\opmod_{0,2,5}$ corresponding to three discs glued together}
\end{figure}
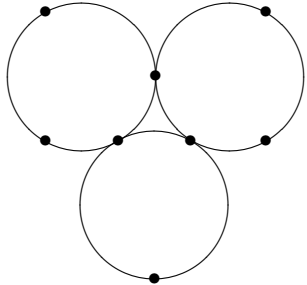

Moduli  spaces $\opmod_{g,n,r}$ of Riemann surfaces with boundary,
and with marked points and possibly nodes on the boundary, are
described. These moduli spaces are
 manifolds with corners; the boundary $\partial \opmod_{g,n,r}$ is the locus of
singular surfaces. These moduli spaces are close relatives of
those studied by Liu and Katz \cite{liu2002, kat_liu2003} and
Fukaya-Oh-Ohta-Ono \cite{fuk_oh_oht_ono2000}; they are an open
subset of the natural Deligne-Mumford compactification of the
moduli of Riemann surfaces with boundary, and are closely related
to the real points of the usual Deligne-Mumford spaces.

The topological type of a curve $C$ in the Deligne-Mumford moduli
space can be described by a graph. In a similar way, the
topological type of a Riemann surface with nodal boundary $\Sigma
\in \opmod_{g,n,r}$ can be described by a type of graph.   If all
of the irreducible components of $\Sigma \in \opmod_{g,n,r}$ are
discs, then the topological type of $\Sigma$ is described by a
ribbon graph, with $r$ external edges. Let $D_{g,n,r} \into
\opmod_{g,n,r}$ be the locus of such surfaces.  I show that the
inclusion $D_{g,n,r} \into \opmod_{g,n,r}$ is a homotopy
equivalence. $\opmod_{g,n,r}$ is a manifold with corners, and so
is homotopy equivalent to its interior $\mc N_{g,n,r}$. When $r =
0$, $\mc N_{g,n,0}$ is homotopy equivalent to $\mc M_{g,n}$, the
space of smooth complex algebraic curves with $n$ marked points.
Therefore there is a homotopy equivalence $D_{g,n,0} \simeq \mc
M_{g,n}$. The locus in $D_{g,n,0}$ of surfaces whose topological
type is described by a fixed ribbon graph is an orbi-cell.  This
immediately implies that the complex of singular chains $C_\ast
(\mc M_{g,n}\otimes \Q)$ is quasi-isomorphic to a complex built
from ribbon graphs, recovering the results of
Harer-Mumford-Thurston-Penner.

The spaces $\opmod_{g,n,r}$ form a modular operad, by gluing
surfaces along marked points.   The statement that
$\opmod_{g,n,r}$ is homotopy equivalent to $D_{g,n,r}$ can be
interpreted as saying that this operad is generated (up to
homotopy) by the moduli spaces $\opmod_{0,1,r}$ of discs. All the
relations also come from $\opmod_{0,1,r}$.  The moduli space
$\opmod_{0,1,r}$ of discs form a cyclic operad -  if two  discs
are glued together, the result is still a (singular) disc.  This
cyclic operad is isomorphic to the topological cyclic operad
$A_\infty^{top}$ of Stasheff \cite{sta1963}.

To make this statement precise I need to develop some operadic
formalism. The notions of cyclic and modular operads, which are
generalizations of operad, were introduced by Getzler and Kapranov
\cite{get_kap1995,get_kap1998}, following Kontsevich's work on
graph cohomology \cite{kon1993,kon1994}. A good introduction to
these operadic concepts can be found in the book of Markl, Shnider
and Stasheff \cite{mar_shn_sta2002}. Roughly, a cyclic (resp.
modular) operad is to a forest (resp. graph) what an operad is to
a rooted forest (recall a forest is a disjoint union of trees).
This can be made precise: I define symmetric monoidal categories
$\RTrees$, $\Trees$ and $\Graphs$, whose morphisms are given by
rooted forests, forests, and graphs. If ${\bf C}$ is a symmetric
monoidal category, a tensor functor $\RTrees \to {\bf C}$ (resp.
$\Trees \to {\bf C}$, $\Graphs \to {\bf C}$) is the same as an
operad (resp. cyclic operad, modular operad) in ${\bf C}$. There
is a functor $\Trees \to \Graphs$; therefore every modular operad
is also a cyclic operad. Conversely, for every cyclic operad $P$,
I define a modular operad $\Mod(P)$, the ``modular envelope'' of
$P$ (by analogy with the universal enveloping algebra of a Lie
algebra). $\Mod(P)$ comes equipped with a map of cyclic operads $P
\to \Mod(P)$. This map is universal, in the sense that if $Q$ is a
modular operad and $P \to Q$ is a map of cyclic operads, there is
a unique map of modular operads $\Mod(P) \to Q$ making the obvious
diagram commute.

The moduli spaces $\opmod_{g,n,r}$ form a topological modular operad $\open$,
by gluing surfaces at marked points.   The moduli spaces $\opmod_{0,1,r}$ form
a cyclic operad, by gluing two discs along marked points to get a singular
disc.   It is known that this cyclic operad is isomorphic to the topological
operad $A_\infty^{top}$ of Stasheff \cite{sta1963}.  Thus there is a map of
topological cyclic operads, $ A_\infty^{top} \to \open$.  By the universal
property of the modular envelope construction, this induces a map $\Mod(
A_\infty^{top}) \to \open$. The main result of this paper then says
\begin{theorem}
The map $\Mod(A_\infty^{top}) \to \open$ is a homotopy equivalence of
orbispace modular operads.
\end{theorem}

Let $C_\ast$ be an appropriate chain-complex functor from topological spaces to
dg $\Q$--vector spaces. Applying $C_\ast$ turns a topological operad into a dg
operad; this gives a dg modular operad $C_\ast(\open)$ and a dg cyclic operad
$C_\ast(A_\infty^{top})$. It is known that $C_\ast(A_\infty^{top})$ is
quasi-isomorphic to the usual algebraic $A_\infty$ operad, $A_\infty^{alg}$.
\begin{corollary}
There is a quasi-isomorphism of dg modular operads
$$
\Mod(A_\infty^{alg} \otimes \Q)  \to C_\ast(\open)
$$
\end{corollary}

The operad $\open$ is the operad controlling open topological
conformal field theory (TCFT); an algebra over the modular operad
$C_\ast(\open)$ is called an  open TCFT. \footnote{This is not the
most general definition of open TCFT; if we generalise to allow
more than one D-brane, we find an open TCFT is the same as an
$A_\infty$ category of Calabi-Yau type. } The reason for this
terminology is as follows. Recall that the operad structure on
$\open$ is given by gluing marked points $P_1 \in
\partial \Sigma_1$, $P_2 \in
\partial \Sigma_2$ together. Replace each point $P_i$ by a small
interval $I_i$; this gives us a homotopy equivalent moduli space.
The operad structure is given by gluing these intervals together.
This is like the string theory operation of gluing an incoming
open string to an outgoing open string.

Applying $\pi_0$ to the homotopy equivalence $\Mod(A_\infty^{top}) \simeq
\open$ recovers the result of Moore and Segal relating topological open field
theory to Frobenius algebras. This is because $\pi_0 (\Mod(A_\infty^{top}))$ is
the modular envelope of the associative cyclic operad; whereas $\pi_0(\open)$
is a modular operad in the category of sets, constructed from isomorphism
classes of topological surfaces with marked intervals on the boundary.  The
result of Moore and Segal can be interpreted as saying that algebras over the
modular operad $\pi_0(\open)$ are precisely associative algebras with inner
product, or in other terms that the operad $\pi_0(\open)$ is isomorphic to
$\Mod( \op{Assoc})$.

The main result of this note is therefore a close analogue of this
theorem of Moore and Segal.  In fact, what is proved here is a
``derived'' version of their result. As Ezra Getzler pointed out
to me, one can interpret the modular operad $\Mod(A_\infty^{alg})$
as being $\mbb {L} \Mod( \op{Assoc})$, where $\mbb {L} \Mod$ is
the left derived functor in the sense of homotopical algebra. This
is because the operad $A_\infty^{alg}$ is a free resolution (and
therefore a cofibrant model) of the operad $\op{Assoc}$ of
associative algebras. Therefore, what is shown here is that $\mbb
{L} \Mod(\op{Assoc} \otimes \Q) \iso C_\ast (\open)$.  One would
hope that in a similar way, the topological $A_\infty$ operad is a
cofibrant model (whatever that means) of the associative cyclic
operad, in the category of topological cyclic operads; this would
show that $\mbb L \Mod (\op{Assoc}) \iso \open$. Of course there
are considerable difficulties making sense of this in the
topological setting.

$A_\infty^{alg} \otimes \Q$ is not just a cofibrant model of the
associative operad over $\Q$; it is the minimal model in the sense
of Markl \cite{mar1996}.   Markl's theory should generalise
without much difficulty to  modular and cyclic operads.  Then
$\Mod(A_\infty^{alg}\otimes \Q)$ is the minimal model of
$C_\ast(\open) \otimes \Q$. As, $\Mod(A_\infty^{alg}\otimes \Q)$
is free (after forgetting the differential), and the image of the
differential is in the space of decomposable elements.  If we
define a homotopy action of a connected modular operad on a
complex to be an action of its minimal model, this should show
that a homotopy action of $C_\ast(\open) \otimes \Q$ on a complex
$V$ with inner product is  the same as an action of
$\Mod(A_\infty^{alg} \otimes \Q)$ on $V$. However, this is the
same as an $A_\infty$ structure on $V$ - that is a homotopy
associative structure. Therefore, we see that an open TCFT
 (at all genera) is precisely a
homotopy associative Frobenius algebra.

I should remark that the operads $\opmod_{g,n,r}$ seem to be
related to the arc operads studied by Penner  and Kauffmann,
Livernet, Penner in the interesting papers \cite{pen2003,
kau_liv_pen2003}. However, they are mostly concerned with the
operads of closed rather than open strings.

\subsection*{Acknowledgements}
I would like to thank Ezra Getzler and Jim Stasheff for their
comments on an earlier version of this note, and  Graeme Segal for
patiently explaining his work to me.

\section{Cyclic and modular operads}
In this section, I give a definition of operads, cyclic operads,
and modular operads.  The definitions presented here have a
slightly different flavour to the usual definitions (although they
are essentially equivalent).  In this note, an operad (resp.
cyclic operad, modular operad) in a symmetric monoidal category
${\bf C}$ is defined to be a tensor functor from a symmetric
monoidal category constructed from rooted forests (resp. forests,
graphs) to ${\bf C}$.   The advantage of this definition is that
it makes the construction of free operads, and of the modular
envelope of a cyclic operad, easier. I begin by defining the
categories $\RTrees$, $\Trees$, $\Graphs$.

A graph is what you think it is :  it is a collection of vertices joined by edges. Graphs may
 be disconnected, there may be external edges (or tails), and vertices may have loops.  Slightly
more degenerate graphs are allowed than is usual: for example, there is a graph
with one vertex and no edges. Here is a picture of a graph. 
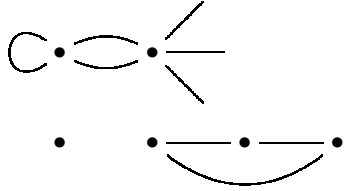
\vspace{5mm}
$$
\xymatrix {
\bullet \ar@(ul,dl)@{-}[]  \ar@/^/@{-}[r] \ar@/_/@{-}[r] & \bullet  \ar@(ur,ur)@{-}[] \ar@(r,r)@{-}[]
\ar@(dr,dr)@{-}[] &  \\
\bullet & \bullet \ar@{-}[r] \ar@/_1.3pc/@{-}[rr]  & \bullet \ar@{-}[r] & \bullet }
$$
\vspace{5mm}

There are various finite sets associated to a graph $\gamma$. There is the set
$C(\gamma)$ of connected components, the set $T(\gamma)$ of tails or external edges,
the set $V(\gamma)$ of vertices, and the set $H(\gamma)$ of germs of edges, or
half-edges.  A half-edge is an edge (internal or external) together with the
choice of a vertex attached to it. There are maps $T(\gamma) \to C(\gamma)$ and
$H(\gamma) \to V(\gamma)$. For a vertex $v$, write $H(v)$ for the fibre of
$H(\gamma) \to V(\gamma)$ at $v$.

For example, in the picture above, $\#C(\gamma) = 3$, $\#T(\gamma) = 3$, $\#V(\gamma) = 6$
and $\#H(\gamma)= 15$.  If $v$ is  the vertex on the upper left, then $\#H(v) = 4$.

It is convenient to give a formal definition of a graph.
\begin{definition}
A \emph{graph} $\gamma$ consists of
\begin{itemize}
\item
A finite set $V(\gamma)$ of vertices;
\item
a finite
set $H(\gamma)$ of half-edges (or germs of edges);
\item
a map $\pi : H(\gamma) \to V(\gamma)$;
\item
an involution $\sigma : H(\gamma) \to H(\gamma)$,
satisfying $\sigma^2 = 1$.
\end{itemize}
The \emph{edges} $E(\gamma)$ of $\gamma$ are the free $\sigma$-orbits in $H(\gamma)$.
The \emph{tails} $T(\gamma)$ of $\gamma$ are the $\sigma$ fixed points of $H(\gamma)$.
The set of connected components $C(\gamma)$ of $\gamma$ is the quotient of
$V(\gamma)$ by the equivalence relation generated by : $v \sim v'$ if there is
a half-edge $h \in H(\gamma)$, with $\pi(h) = v$ and $\pi(\sigma(h)) = v'$.  There is a map
$T(\gamma) \to C(\gamma)$.

The \emph{geometric realization} $\abs \gamma$ of $\gamma$ is the cell complex,
with a $0$-cell for each vertex $v \in V(\gamma)$, and a copy $I_h$ the
interval $I = [0,1]$ for each half-edge $h \in H(\gamma)$.  $0 \in I_h$ is
glued to the vertex $\pi(h)$,  $I_h$ is identified with  $I_{\sigma(h)}$ via
the map $I_h \to I_{\sigma(h)}$, $t \to 1-t$.

A \emph{forest} is a graph all of whose connected components are
contractible, and each of whose vertices is at least trivalent.  A
\emph{rooted forest} is a forest together with a choice of tail
for each connected component, that is a section $C(\gamma) \to
T(\gamma)$ of the projection $T(\gamma) \to C(\gamma)$.

\end{definition}
Now I can define the categories $\Graphs$, $\Trees$ and $\RTrees$.  An object
of $\Graphs$ is a pair $I, J$ of finite sets, together with a map $I \to J$. In
order to distinguish these maps of finite sets from the morphisms in the
category $\Graphs$, I will write $[ I \rrar J ]$ for this object of $\Graphs$.

The morphisms of $\Graphs$ are given by graphs.  A graph $\gamma$ is a morphism
$$
\gamma : [H(\gamma) \rrar V(\gamma)] \to [T(\gamma) \rrar C(\gamma)].
$$
Let $\gamma_1, \gamma_2$ be graphs, with an
isomorphism
$$[T(\gamma_2) \rrar C(\gamma_2)] \iso [H(\gamma_1) \rrar V(\gamma_1)].$$
 Thus, $\gamma_2$ is a morphism
$$\gamma_2 : [H(\gamma_2) \rrar V(\gamma_2)] \to [T(\gamma_2) \rrar C(\gamma_2)],$$
and $\gamma_1$ is a morphism
$$\gamma_1 : [T(\gamma_2) \rrar V(\gamma_2)] \to [T(\gamma_1) \rrar C(\gamma_1)].$$
The composition $\gamma_1 \circ \gamma_2$ can be formed by inserting $\gamma_2$
into $\gamma_1$. That is, each vertex $v \in V(\gamma_1)$ is replaced by the
corresponding connected component of $\gamma_2$ under the identification
$V(\gamma_1) \iso C(\gamma_2)$, and the half-edges $H(v)$ are glued at $v$ to
the corresponding tails of the connected component of $\gamma_2$, using the
identification $H(\gamma_1) \iso T(\gamma_2)$.

The composition of two graphs can be drawn as follows.  Let $\gamma_1$ be the
graph
$$
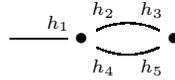
\xymatrix{
\ar@{-}@(l,l)[]^>>{h_1} \ar@{-}@/^/[r]^<<{h_2}^>>{h_3} \ar@{-}@/_/[r]_<<{h_4}_>>{h_5}  \bullet & \bullet
}
$$
and let $\gamma_2$ be the graph
$$
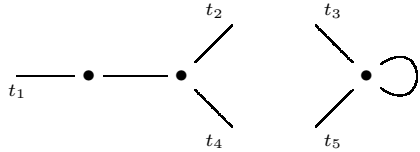
\xymatrix{
\bullet \ar@{-}@(l,l)^{t_1} \ar@{-}[r] & \bullet  \ar@{-}@(ur,ur)^{t_2} \ar@{-}@(dr,dr)_{t_4}
& \hspace{5pt} & \bullet \ar@{-}@(ul,ul)_{t_3} \ar@{-}@(dl,dl)^{t_5} \ar@{-}@(ur,dr) }
$$
In order to form the composition $\gamma_1 \circ \gamma_2$, we need to identify the tails $T(\gamma_2)$
with the half-edges $H(\gamma_1)$.  The tail $t_i \in T(\gamma_2)$ is identified with the half-edge
$h_i \in H(\gamma_1)$.
The connected components of $\gamma_2$ are identified with the
vertices of $\gamma_1$ in the only possible way.  Then $\gamma_1 \circ \gamma_2$
is
$$
\xymatrix{
\bullet \ar@{-}@(l,l) \ar@{-}[r] & \bullet  \ar@{-}@/^/[r] \ar@{-}@/_/[r]
&  \bullet  \ar@{-}@(ur,dr) }
$$

If $[I \rrar J]$ is an object of $\Graphs$, the identity map $[I \rrar J]$ is
given by the graph $\gamma$ with $T(\gamma) = H(\gamma) = I$ and $V(\gamma) =
C(\gamma) = J$. The involution $H(\gamma) \to H(\gamma)$ is the identity, or in
other terms $\gamma$ has no internal edges.     If  $[I \rrar J]$, $[K \rrar
L]$ are two objects of $\Graphs$, an isomorphism (of pairs of finite sets with
morphisms) $[I \rrar J] \iso [K \rrar L]$ corresponds to the graph $\gamma$
with $H(\gamma) = I$, $V(\gamma) = J$, and no internal edges as before; but
with the isomorphism $T(\gamma) \iso K$, $C(\gamma) \iso L$.

$\Graphs$ is a symmetric monoidal category.  The tensor product is defined on objects by
$$
[I \rrar J] \otimes [K \rrar L] = \left[ (I \coprod K) \rrar (J \coprod L) \right]
$$
and on morphisms by
$$\gamma_1 \otimes \gamma_2 = \gamma_1 \coprod \gamma_2.$$

Let $\Trees \subset \Graphs$ be the subcategory such that $\op{Ob}
\Trees \subset \op{Ob} \Graphs$ consists of those $[I \rrar J] \in
\op{Ob} \Graphs$ such that the fibres of the  map of finite sets
$I \rrar J$ are of cardinality at least $3$. $\op{Mor}\Trees
\subset \op{Mor}\Graphs$ consists of forests. $\Trees$ is a
symmetric monoidal category.

Let $\RTrees$ be the category whose objects are maps $I \rrar J$
of finite sets such that the cardinality of the fibres is at least
$3$, together with a section $\sigma : J \to I$.  The morphisms in
$\RTrees$ are given by rooted forests.  The composition $\RTrees$
is defined in a similar way to that in $\Graphs$.  There is a
functor $\RTrees \to \Trees$.

After these preliminaries, I can define the types of operad I need.
\begin{definition}
Let ${\bf C}$ be a symmetric monoidal category. A \emph{modular operad} is a tensor functor
$$P
: \Graphs \to {\bf C} . $$
A \emph{cyclic operad} is a tensor functor $\Trees \to {\bf C}$.  An
\emph{operad} is a tensor functor
$\RTrees \to \bf C$.
\end{definition}
One can see that these definitions are equivalent to the more
traditional definitions. For example, any tensor functor $F :
\RTrees \to \bf C$ is determined by its action on the objects $[(I
\cup \{\ast\}) \rrar \{\ast\}]$, and on the morphisms
corresponding to connected rooted forests with one internal edge.
This is because these objects and morphisms generate $\RTrees$ as
a symmetric monoidal category. For each finite set $I$, let
$$P(I) = F([(I \cup\{\ast\}) \rrar \{\ast\}]) \in \op{Ob} \bf C$$
$P(I)$ is acted on by $\op{Aut} I$.  For each pair $I,J$ of finite
sets, with elements $i \in I$, there is a rooted forest in
$\op{Mor} \RTrees$ with one internal edge, which gives a morphism
$$
[(I \cup \{\ast\}) \rrar \{\ast\}] \otimes [(J \cup \{\ast\}) \rrar \{\ast\}] \to [(I \setminus \{i\}  \cup J
\cup \{\ast\}) \rrar \{\ast\}]
$$
This shows that a tensor functor $\RTrees \to \bf C$ is the same as a collection
of objects $P(I) \in \bf C$, one for each finite set $I$ with at least $3$ elements, together
with an $\op{Aut} (I)$ action on $P(I)$; and for each finite set $J$, and element $i \in I$,
a composition map
$$
P(I) \otimes P(J) \to P( I\setminus \{i\} \cup J )
$$
This composition map must satisfy a certain associativity property; and we find that a tensor
functor $\RTrees \to \bf C$ is the same as an operad in $\bf C$.

Similar remarks hold for cyclic and modular operads. For any
functor $P : \Graphs \to \bf C$ or $\Trees \to \bf C$, and any
finite set $I$, we will abuse notation and write
$$
P(I) = P([I \rrar \ast])
$$
If $n \in \Z_{\ge 0}$, we will also write $P(n)$ for $P([n])$
where $[n] = \{1,2, \ldots, n\}$. For $i \in I, j \in J$, we will
write
$$
\circ_{i,j} : P(I) \otimes P(J) \to P(I \cup J \setminus \{i,j\})
$$
for the composition map coming from the morphism in $\Trees$ or
$\Graphs$,
$$[I \cup J \rrar \{ 1,2 \} ] \to [I \cup J \setminus
\{i,j\} \rrar \{ \ast \}]$$ If $P$ is a functor $\Graphs \to \bf
C$, and $i_1, i_2 \in I$ are distinct, we will write
$$
\circ_{i_1,i_2} : P(I) \to P(I \setminus \{i_1, i_2 \})
$$
for the map coming from the morphism in $\Graphs$,
$$
[I \rrar \{\ast\}] \to [I \setminus \{i_1, i_2 \} \rrar \{\ast\}]
$$
It is easy to see that the composition maps $\circ_{i,j}$ and
$\circ_{i_1,i_2}$ satisfy the usual associativity and equivariance
axioms for a cyclic or modular operad. Also $P : \Graphs \to \bf
C$ or $\Trees \to \bf C$ is completely determined by the $P(I)$,
with their natural $\op{Aut}(I)$ action, and these composition
maps.
\subsection{Free operads and the modular envelope}
Suppose  ${\bf C}_1$, ${\bf C}_2$, ${\bf C}_3$, are categories, and $F : {\bf
C}_1 \to {\bf C}_2$ is a functor. There is a pull-back functor
$$
F^\ast : \op{Fun} ( {\bf C_2}, {\bf C_3} ) \to \op{Fun} ( {\bf C_1},
{\bf C_3})
$$
given by composition with $F$.  In certain nice situations, this
functor admits a left adjoint
$$F_\ast  :\op{Fun} ( {\bf C_1}, {\bf C_3} ) \to \op{Fun} ( {\bf C_2},
{\bf C_3})$$
The construction of $F_\ast$ is an analogue of the familiar induction of group
representations.  If $G : {\bf C_1} \to {\bf C_3}$ is a functor,
then $F_\ast G$  satisfies a universal
property.   There is a morphism of functors
$$
G \to F^\ast F_\ast (G)  = F_\ast(G) \circ F
$$
such that if
$H : {\bf C}_2 \to {\bf C}_3$ is any functor, and
$G \to F^\ast H  $ is a morphism of functors, then there is a unique
morphism of functors $F_\ast G \to H$ such that the diagram
$$
\xymatrix{
G \ar[r] \ar[dr]  & F^\ast F_\ast G \ar[d] \\
& F^\ast H
}
$$
commutes.

A functor $F_\ast G$ with this universal property doesn't always exist;  it
only exists when ${\bf C}_3$ admits enough coproducts. In this case, it is
defined as follows.  For an object $x \in {\bf C}_2$, $F_\ast G(x)$ has a copy
$G(y)_f$ of $G(y)$ for each object $y \in \op{Ob} {\bf C}_1$ with a morphism $f
: F(y) \to x$.  If  $g : y' \to y$ is a morphism in ${\bf C_1}$, then the copy
$G(y')_{f \circ F(g)}$  of $G(y')$ is identified with the image of $G(g) :
G(y') \to G(y) = G(y)_f$.  $F_\ast$, when it exists, is the left adjoint to the
functor $ F^\ast : \op{Fun} (C_2, C_3) \to \op{Fun}(C_1, C_3)$ given by
composition with $F$.

This construction will be applied to construct free modular and cyclic operads,
and to construct a modular operad from a cyclic operad.

Let ${\bf C}$ be a symmetric monoidal category, which is now assumed to be
$k$-linear.  Let $\Pairs$ be the tensor category whose objects are the same as
those of $\Graphs$, that is pairs $I,J$ of finite sets with a map $I \rrar J$.
The morphisms in $\Pairs$ are simply isomorphisms $[I \rrar J] \iso [I' \rrar
J']$, so that $\Pairs$ is a groupoid.  The tensor structure on $\Pairs$ is
given by disjoint union, as before.  Let $P : \Pairs \to {\bf C}$ be a tensor
functor. There is a functor $F : \Pairs \to \Graphs$; in fact $\Pairs$ is the
subcategory of $\Graphs$ whose objects are all those of $\Graphs$ and whose
morphisms are the isomorphisms in $\Graphs$.   $F_\ast P$ is a functor $\Graphs
\to {\bf C}$, that is a modular operad.  $F_\ast P$ is the free modular operad
generated by $P$.

Now suppose $P : \Trees \to {\bf C}$ is a cyclic operad.  Since their is a
functor $\Gamma : \Trees \to \Graphs$, every modular operad is in particular a
cyclic operad. That is, there is a forgetful functor $\Gamma^\ast$ from modular to cyclic
operads. Denote this functor by $\Cyc$. There is a left adjoint $\Gamma_\ast$ to this functor, as long as $\bf C$ admits
enough coproducts.    Denote this functor by $\Mod$.

I call
$\Mod(P)$ the \emph{modular envelope} of $P$, by analogy with the
universal enveloping algebra of a Lie algebra.  The modular envelope satisfies
a universal property. There is a map $P \to \Cyc (\Mod (P))$ of cyclic
operads.  For any modular operad $Q$ with a map $P \to \Cyc (Q)$, there is a
unique map $\Mod (P) \to Q$ such that the diagram
$$
\xymatrix{
P \ar[r] \ar[dr]  & \Cyc (\Mod (P)) \ar[d] \\
& \Cyc (Q)
}
$$
commutes.

Note that there is no problem in extending the  definition of cyclic and
modular operads to the $2$-category of topological orbispaces. Let $\Orb$ be
this $2$ category; there is an obvious functor $F : \Top \to \Orb$. Then, if $P
: \Trees \to \Top$ is a cyclic operad in the category of topological spaces, $F
\circ P: \Trees \to \Orb$ is a cyclic operad in $\Orb$.    Form the modular
envelops $\Mod (P) : \Graphs \to \Top$, the modular envelope of $P$, and
$\Mod (F \circ P)$, the modular envelope of $F \circ P$.  These are not the
same, that is $F \circ \Mod (P) \neq \Mod (F \circ P)$.  In many ways it
is better to consider $\Mod (F \circ P)$, because the construction of
$\Mod$ involves forming quotients by actions of finite groups.

\subsection{Examples of cyclic operads}
The first example is $\mc A$, the associative cyclic operad; this is a cyclic
operad in the tensor category $\Fin$ of finite sets, whose morphisms are
isomorphisms of finite sets.  The tensor structure on $\Fin$ is given by
Cartesian product.

For an object $[I \rrar J] \in \op{Ob} \Trees$, define
$$
\mc A ([I \rrar J]) = \{ \text {cyclic orders on the fibres of } I \rrar J \}
$$
For a morphism $\gamma : [H(\gamma) \rrar V(\gamma)] \to [T(\gamma) \rrar
C(\gamma)]$, we need to define
$$
\mc A(\gamma) : \mc A ([H(\gamma) \rrar V(\gamma)]) \to \mc A ([T(\gamma) \rrar
V(\gamma)])
$$
Note that an element $a \in \mc A( [H(\gamma) \rrar V(\gamma)])$
corresponds to a cyclic order on the set of edges emanating from
each vertex of the forest $\gamma$. That is, $\gamma$ becomes a
ribbon graph.   $\gamma$ can be thickened to form a compact
oriented surface $\Sigma$ with  boundary; as $\gamma$ is a forest,
$\Sigma$ is a disjoint union of discs.  The orientation on
$\partial \Sigma$ induced by that on $\Sigma$ induces a cyclic
order on the tails $T(c)$ for each connected component $c \in
C(\gamma)$, that is a cyclic order on the fibres of $T(\gamma) \to
C(\gamma)$; this defines $\mc A (\gamma)$.

$\mc A$ is called the associative cyclic operad. There is a tensor functor
$\Fin \to \Vect_k$ which sends a finite set $S$ to the vector space $k^{\oplus
S}$ with basis $S$.  The cyclic operad $\mc A$ in $\Fin$ pushes forward to the
usual cyclic operad of associative algebras in the category of vector spaces,
which I also denote by $\mc A$.

The second basic example is $\mc C$, the cyclic operad of commutative algebras.
Again, $\mc C$ is a cyclic operad in the tensor category $\Fin$. For an object
$a \in \op{Ob} \Trees$,  define $\mc C(a) = \{ \ast \}$, the set with one
element. The definition of $\mc C$ on morphisms in $\Trees$ is trivial.

In this paper, the cyclic operad we will be most concerned with is $A_\infty^{alg}$, the operad
of $A_\infty$ algebras.  $A_\infty^{alg}$ is a cyclic operad in the category of differential
graded $\Q$ vector spaces, with differential of degree $-1$. (This choice of degree of the differential
is so as to be consistent with the choice for chain complexes of topological spaces later on).
 Let $\dgv_\Q$ be this tensor category, and let $\Vect_\Q^{\Z}$
be the category of graded $\Q$ vector spaces. As a graded cyclic operad, that is forgetting the differential, $A_\infty^{alg}$
is freely generated by a certain tensor functor $F : \op{Fin}  \to \Vect_\Q^\Z$.
On the finite set $[I \rrar \ast]$, $F$ is defined to be the $\Q$ vector
space with basis the set of cyclic orders on $I$, situated in degree $3 - \#I$.  The functor
$F$ is extended to $\Fin$ by making it a tensor functor.

Taking the free cyclic operad on $F$ gives us a tensor functor
$A_{\infty \#}^{alg} : \Trees \to \Vect_\Q^\Z$.  The cyclic operad
$A_\infty^{alg}$ is obtained from this by adjoining a certain
differential. If $[I \rrar J] \in \op{Ob} \Trees$, then
$A^{alg}_{\infty \#}([I \rrar J])$  has a basis corresponding to
forests $\gamma$ with isomorphisms $[T(\gamma) \rrar C(\gamma)]
\iso [I \rrar J]$,  with cyclic orders on the fibres of
$[H(\gamma) \rrar V(\gamma)]$, and with orderings of the set
$V(\gamma)$. Reordering the set $V(\gamma)$ changes the basis
element by a sign corresponding to the order of the permutation.
The basis element corresponding to a forest $\gamma$ has degree
$\#H(\gamma) - 3\#V(\gamma)$.  The differential is defined on
these basis elements by summing over all ways to add an edge to
the forest, with appropriate sign.

\subsection{Algebras over cyclic and modular operads}
Let $P$ be a cyclic or modular operad in the tensor category ${\bf C}$; ${\bf
C}$ is now assumed to be one of the categories of finite dimensional vector
spaces, $\Z$-graded vector spaces or dg vector spaces.

I want to define the notion of a $P$-action on an object $V \in
\op{Ob} {\bf C}$.  Suppose $\diag \in V^{\otimes 2}$ is a closed
element of degree $0$ say. Then I define a modular operad
$\op{End} (V, \diag)$, as follows. For an object $[I \rrar J] \in
\op{Ob} \Graphs$, define
$$
\op{End}(V , \diag) ([I \rrar J]) = V^{\vee \otimes I}
$$
The modular operad structure on $\op{End}(V, \diag)$ uses the
tensor $\diag$. If $\gamma : [H(\gamma) \rrar V(\gamma)] \to
[T(\gamma) \rrar C(\gamma)]$ is a morphism in $\Graphs$, then
define
$$
\op{End}(V, \diag) (\gamma) : V^{\vee \otimes H(\gamma)}  \to
V^{\vee \otimes T(\gamma)}
$$
to be
$$
\op{End}(V, \diag) (\gamma) = \otimes_{e \in E(\gamma)} \diag_e
$$
that is the tensor product of a copy of $\diag$ for each edge $e
\in E(\gamma)$, acting on the half-edges corresponding to $e$.
With these definitions, it is easy to see that the composition
maps are of degree $0$.

If $P$ is a cyclic operad, then a $P$ action on  $(V, \diag)$ is a
map of cyclic operads
$$
P \to \Cyc \left( \op{End} (V, \diag) \right)
$$
If $P$ is a modular operad, a $P$ action on $(V, \diag)$ is a map
of modular operads
$$
P \to \op{End} (V, \diag)
$$
Note that if $P$ is a cyclic operad, then a $P$ action on $(V,
\diag)$ is the same as a $\Mod (P)$ action on $(V, \diag)$.

The definition of action of a cyclic operad on a complex presented
here is possibly too restrictive; for interesting work on
generalising the notion of an algebra over a cyclic operad see
\cite{tra2001, lon_tra2003}.

\section{The open TCFT operad and the $A_\infty$ operad}

\subsection{Recollections on Riemann surfaces with boundary}
A Riemann surface of genus $g$ with $n> 0$  boundary components has the following
equivalent descriptions.
\begin{enumerate}
\item
A compact connected ringed space $\Sigma$, isomorphic as a topological space to
a genus $g$ surface with $n$ boundary components, and locally isomorphic to $\{
z \in \C \mid \op{Im} z \ge 0\}$, with its sheaf of holomorphic functions.
\item
A smooth, proper, connected, complex algebraic curve $C$ of genus $2g-1+n$,
with a real structure, such that $C \setminus C(\R)$ has precisely two
components, and $C(\R)$ consists of $n$ disjoint circles; together with a
choice of a component of $C \setminus C(\R)$.
\item
Suppose $2g - 2 + n > 0$. Then, a Riemann surface with boundary is equivalently
a 2-dimensional connected compact oriented $C^\infty$ manifold $\Sigma$ with
boundary, of genus $g$ with $n$ boundary components, together with a metric of
constant curvature $-1$ such that the boundary is geodesic.
\end{enumerate}
(2) and (3) can be shown to be equivalent (when $2g - 2 + n > 0$) as follows.
Given $\Sigma$,  $C$ is obtained by gluing $\Sigma$ and $\overline \Sigma$
along their boundary.  Conversely, given $C$, $\Sigma$ is the closure of the
chosen component of $C\setminus C(\R)$ in $C$.  The hyperbolic metric on
$\Sigma$ is the restriction of the unique complete hyperbolic metric on $C$
compatible with the complex structure.

I will also need nodal Riemann surfaces with boundary.  To define these it is
easiest to use the algebraic description (2).  A nodal Riemann surface with
boundary is a proper algebraic curve $C$, with at most nodal singularities, and
a real structure. The real structure on each connected component $C_0$ of the
normalization $\widetilde C$ of $C$ must be of the form $(2)$ above; we also
require  a choice of component of $C_0 \setminus C_0(\R)$. All the nodes of $C$
are required to be real, that is in $C(\R)$. Let $\Sigma$ be the closure in $C$
of the chosen components of $C \setminus C(\R)$; $\Sigma$ is a Riemann surface
with nodal boundary.  Near a node, $\Sigma$ looks like
$$
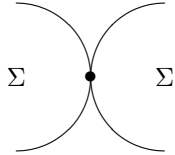
\begin{xy}
0*{\Sigma}, 0*\cir<30pt>{r^l}, 0,a(0), **{}, 0+/28pt/*{\bullet},
0+/56pt/*\cir<30pt>{l^r}, 0+/56pt/*{\Sigma}
\end{xy}
$$

The number of boundary components of $\Sigma$ can be defined as follows.
$\partial \Sigma$ will be a  union of circles, glued
together at points as above. Define a smoothing of $\partial \Sigma$, by replacing each node as
above by
$$
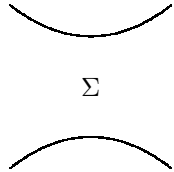
\xymatrix{
\ar@/_1.2pc/@{-}[rr] &  &  \\
 & \Sigma &  \\
\ar@/^1.2pc/@{-}[rr]&  &
}$$
The number of boundary components of $\Sigma$ is defined to be
the number of connected components of this smoothing.

$\Sigma$ has genus $g$ if it has $n$ boundary components and the genus of the nodal algebraic
curve $C = \Sigma_{\cup \partial \sigma} \overline \Sigma$ is $2g-1+n$.

\subsection{The open TCFT operad}

For a finite set $I$ and an integer $n$,  let $\opmod_{g,n,I}$ be the moduli
space of Riemann surfaces $\Sigma$ of genus $g$ with boundary, possibly with
nodes on the boundary, with $n$ boundary components, and with $I$ marked points
on the boundary. The associated algebraic curve $C$, obtained from gluing
$\Sigma$ and $\overline \Sigma$, must be stable, and of genus $2g-1+n$.
Stability is equivalent to the statement that there are only finitely many
automorphisms of $\Sigma$ preserving the marked points.   Let $\mc N_{g,n,I}
\subset \opmod_{g,n,I}$ be the locus of non-singular Riemann surfaces (with
boundary).  These moduli spaces were first constructed by Liu in \cite{liu2002}.
\begin{lemma}
$\opmod_{g,n,I}$ is an orbifold with corners of dimension $6g - 6 + 3n + \#I$.
The interior of $\opmod_{g,n,I}$ is $\mc N_{g,n,I}$. The inclusion $\mc
N_{g,n,I} \into \opmod_{g,n,I}$ is a homotopy equivalence.
\end{lemma}

The spaces $\opmod_{g,n,I}$ form a modular operad $\open$.  For an object $[I
\rrar J] \in \op{Ob} \Graphs$,  define
$$
\open([I\rrar J]) = \prod_{j \in J} \left( \coprod_{g,n}
\opmod_{g,n,I_j} \right)
$$
Let $\gamma : [H(\gamma) \rrar V(\gamma)] \to [T(\gamma \rrar
C(\gamma)]$ be a map. An element in $\open ([H(\gamma) \rrar V(\gamma)])$ corresponds to
a Riemann surface $\Sigma_v$ for each vertex $v \in V(\gamma)$, and a marked point on $\Sigma_v$
for each half-edge at $v$.   Define
$$
\open(\gamma) : \prod_{v \in V(\gamma)} \left( \coprod_{g,n}
\opmod_{g,n,H(v)} \right) \to \prod_{c \in C(\gamma)} \left( \coprod_{g,n}
\opmod_{g,n,T(c)} \right)
$$
by gluing the disconnected surface $\Sigma = \coprod \Sigma_v$ corresponding to
a point in \linebreak $\prod_{v \in V(\gamma)} \left( \coprod_{g,n}
\opmod_{g,n,H(v)} \right)$ to itself, using the edges of $\gamma$ to identify
marked points.

I also need a sub cyclic operad of $\open$, which will be identified with
Stasheff's topological $A_\infty$ operad.  For $[I \rrar
  J] \in \op{Ob} \Trees$, define
$$
A_\infty^{top} ([I \rrar J]) = \prod_{j \in J} \opmod_{0,1,I_j}
$$
This definition extends in an obvious way to a functor $\Trees \to
\Top$, defining a topological cyclic operad $A_\infty^{top}$. The spaces $\opmod_{0,1,I}$
have a natural orientation. As, the open part $\mc N_{0,1,I} \subset \opmod_{0,1,I}$ can be
identified with the quotient of the space of $\#I$ distinct points on the oriented circle $S^1$ by the
action of $\op{PSL}_2(\R)$. The natural orientation on $\op{PSL}_2(\R)$, coming from it's simply
transitive action on the set of unit tangent vectors to the upper half plane, induces an orientation
on $\mc N_{0,1,I}$ and so on $\opmod_{0,1,I}$.

\begin{proposition}
$A_\infty^{top}$ is isomorphic as an operad to the topological $A_\infty$
  operad of Stasheff.
\end{proposition}
\begin{proof}
This result is well-known to experts.  The compactifications of
moduli spaces of marked points on the boundary of a disc used here
are the same as those used in Lagrangian Floer homology.  This
result is therefore the reason for $A_\infty$ relations holding in
the Fukaya category.  A proof is presented in \cite{dev1999}.

I'll briefly sketch a proof. For a cyclically ordered finite set
$I$, let $D_I$ be the compactified moduli space of discs with $I$
marked points on the boundary, such that the natural cyclic order
on the boundary coincides with the given one on $I$. So $D_I$ is a
connected component of $\opmod_{0,1,I}$. Fix three consecutive
elements $0,1,\infty \in I$, which we put at $ 0,1,\infty$ on the
disc. Then, we can identify $D_I$ with a compactification of the
space of $I\setminus \{0,1,\infty\}$ points on the interval.
Further, $D_I$ has a cell decomposition, with cells labelled by
rooted ribbon-trees which are at least tri-valent. The open cells
are given by singular discs of fixed topological type; the root is
given by $1$, say. Now it is not difficult to identify this cell
complex with an associahedron.
\end{proof}

The cyclic operad $A_\infty^{top}$ has a cell decomposition, with
cells labelled by ribbon forests.  The open cells, as before, are
given by the surfaces of fixed topological type. The open cells
have a natural orientation.  This cell decomposition is compatible
with the cyclic operad structure. Let $A_\infty^{alg}$ be the dg
cyclic operad obtained from the cellular chain complexes of
$A_\infty^{top}$.   $A_\infty^{alg}$ is the standard $A_\infty$ dg
cyclic operad, as one can see easily using our earlier description
of the latter in terms of forests. The main point is that the
boundary of $\opmod_{0,1,n}$, with appropriate orientation, is the
sum of copies of $\opmod_{0,1,n_1} \times \opmod_{0,1,n_2}$, where
$n_1 + n_2 - 2 = n$, with appropriate signs.

Let $C_\ast$ be an appropriate chain complex, which has a K\"unneth
map $C_\ast (X) \otimes C_\ast(Y) \to C_\ast(X \times Y)$.  Then
$C_\ast(\open)$ is a dg modular operad.

\subsection{The open TCFT operad and the $A_\infty$ operad}
In this section I show
\begin{theorem}
  There is a homotopy equivalence of orbi-space modular operads
\begin{equation*}
  \open \simeq \Mod(A_\infty^{top})
\end{equation*}
\label{main theorem}
\end{theorem}
$\Mod(A_\infty^{top})$ is to be considered as an orbi-space.  That is,
consider $A_\infty^{top}$ as a cyclic operad in $\Orb$, the $2$ category of
orbi-spaces, and then apply $\Mod$.

An immediate corollary of this theorem is
\begin{corollary}
There is a quasi-isomorphism of dg modular operads over $\Q$,
$$
\Mod ( A_\infty^{alg} \otimes \Q ) \iso C_\ast( \open)
\otimes \Q
$$
\end{corollary}
Markl \cite{mar1996} introduced the notion of \emph{minimal model}
of a dg operad.  One should be able to generalise this definition
to modular operads, and   show that $\Mod(A_\infty^{alg}\otimes
\Q)$ is the minimal model for $C_\ast(\open)\otimes \Q)$. Indeed,
$\Mod(A_\infty^{alg} \otimes \Q)$ is free as a graded modular
operad (forgetting the differential), and the image of the
differential consists of decomposable elements.  Let $(V, \diag)$
be a complex with an element $\diag \in V^{\otimes 2}$. Following
Markl, one could define a homotopy action of the dg modular operad
$C_\ast(\open)$ on $(V,\diag)$ as an action of the minimal model
on $(V,\diag)$. Thus we see, that a action of $C_\ast(\open)$ on
$(V,\diag)$ is the same as a homotopy class of homotopy
associative structure on $(V,\diag)$.

The first step in the proof of theorem \ref{main theorem} is to construct a map
of modular operads $\Mod(A_\infty^{top}) \to \open$. By the universal
property of $\Mod$, it is sufficient to give a map of cyclic operads
$A_\infty^{top} \to \open$; but such a map has already been defined. Let $\Phi
: \Mod (A_\infty^{top}) \to \open$ be the resulting map of modular
operads.

\begin{proposition}
$\Phi$ is a homotopy equivalence.
\label{qis}
\end{proposition}
Let $D_{g,n,I} \subset \opmod_{g,n,I}$ be the locus consisting of curves, with
nodes at the boundary, each of whose irreducible components is a disc.  One can
easily show that $D_{g,n,I} \iso \Mod(A_\infty^{top})$, and that the map
$D_{g,n,I} \to \opmod_{g,n,I}$ is just the map $\Phi$ described above.
\begin{proposition}
The inclusion $D_{g,n,I} \into \opmod_{g,n,I}$ is a homotopy equivalence.
\label{htpy equiv}
\end{proposition}
This implies proposition \ref{qis}.
 To prove that $D_{g,n,I} \into \opmod_{g,n,I}$ is a homotopy equivalence, it is sufficient to show that
\begin{proposition}
For $(g,n) \neq (0,1)$, the inclusion $\partial \opmod_{g,n,I}
\into \opmod_{g,n,I}$ is a homotopy equivalence. \label{htpy equiv
boundary}
\end{proposition}
We will prove this as long as $(g,n) \neq (0,2)$; this case is
easy and is left to the reader.

The idea of the proof is very simple. Given a surface $\Sigma \in
\mc N_{g,n,I}$, we use the canonical hyperbolic metric on $\Sigma$
to flow $\partial \Sigma$ inwards; eventually, $\Sigma$ becomes
singular, and we construct a deformation retract of
$\opmod_{g,n,I}$ onto it's boundary.

 It is easier to apply this
procedure to  non-singular surfaces. Therefore, the first step is
to move the boundary $\partial \opmod_{g,n,I}$ inwards a little
bit.

Let $T$ be a tubular neighbourhood of the boundary $\partial
\opmod_{g,n,I}$ which is locally isomorphic to $\partial
\opmod_{g,n,I} \times [0,1)$.  Let $\mc N'_{g,n,I} =
\opmod_{g,n,I} \setminus T$. $\mc N'_{g,n,I}$ is a manifold with
boundary, and the pair $(\mc N'_{g,n,I}, \partial \mc N'_{g,n,I})$
is homotopy equivalent to the pair $(\opmod_{g,n,I}, \partial
\opmod_{g,n,I})$. Therefore it is sufficient to show that the
inclusion $\partial \mc N'_{g,n,I} \into \mc N'_{g,n,I}$ is a
homotopy equivalence.

Now we describe a map $$\Phi : \mc N'_{g,n,I} \times [0,1] \to
\mc N'_{g,n,I}$$ which is a deformation retraction of the
inclusion $\partial \mc N'_{g,n,I} \into \mc N'_{g,n,I}$.

The map $\Phi$ is constructed as follows. Each surface $\Sigma \in
\mc N'_{g,n,I}$ has a canonical metric of constant curvature $-1$
with geodesic boundary.  Use the exponential map of this metric,
and the inward pointing unit normal to $\partial \Sigma \into
\Sigma$, to flow the boundary $\partial \Sigma$ inwards.  For $t
\in \R_{\ge 0}$, let $\Sigma_t$ be the surface with boundary
obtained by flowing in $\partial \Sigma$ a distance $t$.
Eventually, the boundary intersects itself, and we end up with a
surface in $\partial \opmod_{g,n,I}$; before this happens, we must
have hit $\partial \mc N'_{g,n,I}$. More precisely,
\begin{lemma}
There is a unique $S \in \R_{\ge 0}$ such that $\Sigma_S \in
\partial \mc N'_{g,n,I}$, and $\Sigma_t$ is in the interior of
$\mc N'_{g,n,I}$ for all $t < S$.
\end{lemma}
The map $\Phi$ is now defined by  $$\Phi(\Sigma, x) =
\Sigma_{Sx}$$
\begin{proof}[Proof of lemma]
It is sufficient to show that for some $T$, $\Sigma_T$ (after
forgetting the marked points) is in $\partial \opmod_{g,n,0}$.
This will imply that the family of surfaces will have passed
through $\partial \mc N'_{g,n,I}$.

Let $T$ be the first time at which $\Sigma_T$ is singular. We have
to check that $\Sigma_T \in \partial \opmod_{g,n,0}$ (after
forgetting the marked points).

By doubling $\Sigma$, one can see that for any path $\phi : [0,1]
\to \Sigma$, with $\phi (\{0,1\}) \subset \partial \Sigma$, there
is a unique geodesic $\gamma$, homotopy equivalent to $\phi$
relative to $\partial \Sigma$, which is normal to $\partial
\Sigma$. Further, $\gamma$ minimises length in this homotopy
class.

Let $p_1, p_2 \in \partial \Sigma$ be two distinct points. They
collide at time $t$ if there is a point $x \in \Sigma$, and
geodesics $(p_1, x)$, $(p_2, x)$ of length $t$.  If the piecewise
geodesic $(p_1,x)(x,p_2)$ is not an actual geodesic, then there
must be some geodesic $\gamma$, in the same homotopy class as
$(p_1,x)(x,p_2)$, which is normal to $\partial \Sigma$ and of
length shorter than $2t$.  This will imply that $\Sigma$ will have
become singular before time $t$; therefore, at the first time $T$
at which $\Sigma$ becomes singular, all these piecewise geodesics
are smooth.

The time $T$ is the half the minimum length of a geodesic $\gamma
: (I,
\partial I) \to (\Sigma, \partial \Sigma)$ which is normal to the
boundary.  In order to show that $\Sigma_T \in \partial
\opmod_{g,n,0}$, we have to check 2 things.
\begin{enumerate}
\item
There are no three distinct points $p_1, p_2, p_3, \in \partial
\Sigma$, which collide at time $T$ and at the same point $x \in
\Sigma$. That is, $\Sigma_T$ has precisely nodal singularities.
\item
$\Sigma_T$ is stable : there are no irreducible components of
$\Sigma_T$ which are discs with $\le 2$ nodes.
\end{enumerate}
Note that as we forget the marked points of $\Sigma_T$, it doesn't
matter if they collide with the nodes or each other.

For the first point, suppose $p_1, p_2, p_3 \in
\partial \Sigma$ had this property.  Then, there is are geodesics
$(p_i, x)$  of length $T$, and we have seen that each piecewise
smooth geodesic $(p_i,x)(x_,p_j)$ is necessarily smooth. That is,
the tangent vectors to each $(p_i, x)$ at $x$ all coincide.  This
is impossible, as two geodesics which are tangential at any point
are the same.

The second point is clear : if we split off a disc with $2$ nodes,
then we would find two geodesics $\gamma_1, \gamma_2$, normal to
the boundary, distinct, and in the same homotopy class.  If we
split off a disc with one node, then we would find a contractible
geodesic $\gamma$ of positive length, and normal to the boundary.

\end{proof}

\bibliographystyle{plain}
\bibliography{bibliograph}

\end{document}